\documentclass[12pt]{article}
\usepackage{amssymb,amsthm,amsmath,latexsym}
\usepackage{url}
\usepackage{color}
\usepackage{comment}
\usepackage{lscape}
\usepackage{mathtools}
\usepackage{enumitem}
\usepackage{graphicx}
\newtheorem{thm}{Theorem}
\newtheorem{prop}[thm]{Proposition}

\newtheorem{cor}[thm]{Corollary}
\theoremstyle{remark}

\newcommand{\FF}{\mathbb{F}}

\newcommand{\CC}{\boldsymbol{C}}

\newcommand{\ww}{\omega}
\newcommand{\vv}{\omega^2}

\DeclareMathOperator{\wt}{wt}

\newcommand{\Rep}{\operatorname{Rep}_3}
\allowdisplaybreaks

\begin{document}

\title{Determination of weight enumerators of binary extremal self-dual 
codes of length $60$}

\author{
Masaaki Harada\thanks{
Research Center for Pure and Applied Mathematics,
Graduate School of Information Sciences,
Tohoku University, Sendai 980--8579, Japan.
email: \protect\url{mharada@tohoku.ac.jp}
}
}

\maketitle

\begin{abstract}
We construct three binary extremal self-dual codes of length $60$
with weight enumerators for which no binary extremal self-dual 
codes of length $60$ were previously known to exist.
This completes the determination of the weight enumerators
that occur for binary extremal self-dual codes of length $60$.
\end{abstract}

\section{Introduction}
Self-dual codes are one of the most interesting classes of codes.
This interest is justified by many combinatorial objects
and algebraic objects related to self-dual codes
(see, e.g., \cite{SPLAG} and \cite{RS-Handbook}).

Binary self-dual codes are divided into two classes.
A binary self-dual code $C$ is \emph{doubly even} if the weights of all
codewords of $C$ are multiples of $4$, and \emph{singly even} if there is at 
least one codeword of weight $\equiv 2 \pmod 4$.
For a positive integer $n$, a binary self-dual code of length $n$
exists if and only if $n$ is even, and a doubly even self-dual code
of length $n$ exists if and only if $n \equiv 0 \pmod 8$.
Shadows for binary self-dual codes were introduced by Conway and
Sloane~\cite{C-S} to derive upper bounds on the minimum weight
and restrictions on the weight enumerators of binary singly even
self-dual codes.
They also gave the largest possible minimum weights and the
corresponding possible weight enumerators for lengths up to $64$
and length $72$ (see~\cite[Table~I]{C-S} for the minimum weights
and see also~\cite{GH60} for the weight enumerators at length $60$).
Subsequently, using shadows, Rains~\cite{Rains} improved the general
upper bound on the minimum weight.
He showed that the minimum weight $d$ of a binary self-dual code
of length $n$ satisfies $d\le 4\lfloor n/24\rfloor+4$,
unless $n\equiv22\pmod{24}$, in which case
$d\le 4\lfloor n/24\rfloor+6$.
A binary self-dual code meeting this bound is called \emph{extremal}.

It is a fundamental problem to determine which of the possible
weight enumerators given in~\cite{C-S} actually occur.
In this note, we consider this problem for binary extremal
self-dual codes of length $60$.
The length $60$ is of particular interest for binary singly even self-dual codes with 
minimum weight $12$.
The smallest length $n$ for which a binary singly even self-dual
$[n,n/2,12]$ code exists is either $56$ or $60$~\cite[Table~I]{C-S}.
It is not known whether a binary singly even self-dual $[56,28,12]$ code exists.
Thus, $60$ is currently the smallest length for which binary singly even
self-dual codes with minimum weight $12$ are known to exist
(see Section \ref{Sec:R} for another reason for our interest in length $60$).

The possible weight enumerators of binary extremal
self-dual codes of length $60$ are known as follows:
\begin{eqnarray*}
W_{60,1} &=& 1 +(2555+64 \beta) y^{12} + (33600 - 384 \beta)y^{14} + \cdots,
\\
W_{60,2} &=& 1 + 3451y^{12} + 24128 y^{14} + \cdots,
\end{eqnarray*}
where $\beta$ is an integer. 
A binary extremal self-dual code with weight enumerator $W_{60,2}$
was found in~\cite{C-S}.
If there exists a binary extremal self-dual
$[60,30,12]$ code with weight enumerator $W_{60,1}$,
then 
\[
\beta \in \{0,1,2,\ldots,8,10\}
\]
\cite{HM}.
For $\beta=0,1,2,5,6,7$ and $10$,
a binary extremal self-dual code with weight enumerator $W_{60,1}$
was found in~\cite{TJ98}, \cite{BRY}, \cite{H60}, \cite{YL}, \cite{H60}, 
\cite{DH03} and \cite{GH60}, respectively (see also~\cite{Roberts}).

The main aim of this note is to show the following:

\begin{prop}\label{main}
For each $\beta\in\{3,4,8\}$, there exists a binary extremal 
self-dual $[60,30,12]$ code with weight enumerator $W_{60,1}$.
\end{prop}

As a consequence, we have the following:

\begin{cor}
There exists a binary extremal self-dual $[60,30,12]$
code with weight enumerator $W_{60,1}$
if and only if $\beta \in \{0,1,2,\ldots,8,10\}$.
\end{cor}

This completes the determination of the weight enumerators
that occur for binary extremal self-dual codes of length $60$.

\section{Preliminaries}\label{Sec:2}

Let $\FF_q$ denote the finite field of order $q$, where $q$ is a prime power.
We take $\FF_2$ to be $\{ 0,1\}$, and
we take $\FF_4$ to be $\{ 0,1,\ww, \ww^2\}$, where $\omega^2 = \omega +1$.
An $[n,k]$ \emph{code} $C$ over $\FF_q$
is a $k$-dimensional vector subspace of $\FF_q^n$.
Codes over $\FF_2$ and $\FF_4$ are called 
\emph{binary} and \emph{quaternary}, respectively.
The parameter $n$ is called the \emph{length} of $C$.
The \emph{weight} $\wt(x)$ of a vector $x \in \FF_q^n$ is
the number of nonzero components of $x$.
A vector of $C$ is called a \emph{codeword} of $C$.
The minimum weight among the nonzero codewords of $C$ is called
the \emph{minimum weight} of $C$.
An $[n,k,d]$ code
is an $[n,k]$ code with minimum weight $d$.
The \emph{weight enumerator} of $C$
is given by the polynomial $\sum_{c \in C}y^{\wt(c)}$.
Two $[n,k]$ codes $C$ and $C'$ over $\FF_q$ are
\emph{equivalent}, denoted $C \cong C'$,
if there is an $n \times n$ monomial matrix $P$ over $\FF_q$ with
$C' = CP$, where $CP=\{ x P \mid x \in C\}$.
The \emph{automorphism group} of $C$ consists of all
monomial matrices $P$ with $C = CP$.

The \emph{dual} code $C^{\perp}$ of a binary code $C$ of length $n$
is defined as follows:
\[
C^{\perp}=
\{x \in \FF_2^n \mid \langle x,y\rangle = 0 \text{ for all } y \in C\},
\]
where
$
\langle x,y\rangle = \sum_{i=1}^{n} x_i {y_i}
$
for $x=(x_1,x_2,\ldots,x_n)$ and $y=(y_1,y_2,\ldots,y_n) \in \FF_{2}^n$.
In addition, 
the \emph{Hermitian dual} code $C^{\perp_H}$ of a quaternary code $C$ of length $n$
is defined as follows:
\[
C^{\perp_H}=
\{x \in \FF_{4}^n \mid \langle x,y\rangle_H = 0 \text{ for all } y \in C\},
\]
where
$
\langle x,y\rangle_H= \sum_{i=1}^{n} x_i y_i^2
$
for $x=(x_1,x_2,\ldots,x_n)$ and $y=(y_1,y_2,\ldots,y_n) \in \FF_{4}^n$.
A binary code $C$ is called
\emph{self-dual} if $C=C^\perp$.
A quaternary code $C$ is called
\emph{Hermitian self-dual} if $C=C^{\perp_H}$.

Let $C$ be a binary singly even self-dual code.
Let $C_0$ denote the 
subcode of $C$ consisting of 
codewords having weight $\equiv0\pmod4$.
The \emph{shadow} $S$ of $C$ is defined to be $C_0^\perp \setminus C$.
Shadows for self-dual codes were introduced by Conway and
Sloane~\cite{C-S}.

\section{Extremal self-dual codes of length 60}\label{Sec:Result}

In this section, we describe how to construct a binary extremal self-dual $[60,30,12]$
code with weight enumerator $W_{60,1}$ for $\beta=3,4$ and $8$.

We employ the \emph{cubic construction} described in~\cite{Cubic}, \cite{CKO} and \cite{LS}.
Our formulation differs from the original construction only
in the ordering of coordinates and yields an equivalent code.
Let $C_2$ be a binary self-dual code of length $20$
and let $C_4$ be a quaternary Hermitian self-dual code of the same length.
We define a map $\phi$ from $\FF_4$ to $\FF_2^3$ as follows:
\[
\phi(0)=(0,0,0),
\phi(1)=(0,1,1),
\phi(\ww)=(1,0,1) \text{ and }
\phi(\vv)=(1,1,0).
\]
We also define a map $\Rep$ from $\FF_2^{20}$ to $\FF_2^{60}$ as follows:
\[
\Rep(x)=(x_1,x_1,x_1,x_2,x_2,x_2,\ldots,x_{20},x_{20},x_{20})
\]
for $x=(x_1,x_2,\ldots,x_{20}) \in \FF_{2}^{20}$.
We extend $\phi$ coordinatewise to a map from $\FF_4^{20}$ to $\FF_2^{60}$,
also denoted by $\phi$, by setting
\[
\phi(y)=\bigl(\phi(y_1),\phi(y_2),\ldots,\phi(y_{20})\bigr)
\]
for $y=(y_1,y_2,\ldots,y_{20})\in\FF_4^{20}$.
Then
\begin{equation*}
\CC(C_2,C_4)=\{\Rep(x)+\phi(y)\mid x\in C_2, y\in C_4\}
\end{equation*}
is a binary self-dual code of length $60$.

The classifications of binary self-dual codes of length $20$ and quaternary Hermitian
self-dual codes of length $20$ were given in~\cite{P} and \cite{HM-F4}, respectively.
Namely,
there exist $16$ inequivalent binary self-dual codes of length $20$ and
there exist $3427$ inequivalent quaternary Hermitian self-dual codes of length $20$.
Generator matrices of these codes can be obtained electronically from~\cite{HM-DB}.
We denote by $C_{2,i}$ and $C_{4,j}$
the $i$-th binary self-dual code and the $j$-th
quaternary Hermitian self-dual code in~\cite{HM-DB}, respectively.
As noted in~\cite[p.~2254]{Cubic}, the cubic construction yields numerous binary 
self-dual codes of length $60$ by varying the constituent codes and coordinate 
permutations.
Among these codes, we found extremal self-dual
$[60,30,12]$ codes 
with weight enumerators $W_{60,1}$ for $\beta=3,4$ and $8$.

We give explicit constructions of these codes.
Let $S_{20}$ denote the symmetric group of degree $20$.
For $\sigma\in S_{20}$ and a binary code $C_2$ of length $20$, let $\sigma(C_2)$ denote the code obtained
from $C_2$ by moving coordinate $i$ to position $\sigma(i)$.
Define the following binary self-dual codes of length $60$:
\begin{align*}
C_{2,60,3}&=\CC(\sigma_3(C_{2,14}),C_{4,3168}),\\
C_{2,60,4}&=\CC(\sigma_4(C_{2,15}),C_{4,3168})\text{ and }\\
C_{2,60,8}&=\CC(\sigma_8(C_{2,12}),C_{4,439}),
\end{align*}
where
\begin{align*}
\sigma_3
 &= (1,20,4,14,19,17,15,8,10,18)
(2,13,16,7,6,9,12)(3,11),\\
\sigma_4
 &= (1,9)(2,10,14,19,17,18,7,3,6,16,5,12,4,20)
    (8,15,13,11)\text{ and }\\
\sigma_8
 &= (1,20,9,14,16)(2,8,17)(3,13,5)
    (4,18,7)(6,11,15)(10,12).
\end{align*}
We verified by \textsc{Magma}~\cite{Magma} that
$C_{2,60,3}$, $C_{2,60,4}$ and $C_{2,60,8}$
have minimum weight $12$ and weight enumerators
$W_{60,1}$ for $\beta=3,4$ and $8$, respectively.
By~\cite[Theorem~5.6]{Cubic}, the shadow of $\CC(C_2,C_4)$
is obtained by replacing $C_2$ with its shadow in the cubic construction.
In searching for these codes, we used this description to guide
the choice of constituents and coordinate permutations, seeking
minimum weight $12$ and exactly $\beta$ shadow vectors of weight $6$
simultaneously for $\beta\in\{3,4,8\}$.
We verified by \textsc{Magma}~\cite{Magma} that
the automorphism groups of $C_{2,60,3}$, $C_{2,60,4}$ and
$C_{2,60,8}$ have orders $12$, $6$ and $3$, respectively.

For convenience, we give explicit generator matrices in standard form
for codes equivalent to $C_{2,60,\beta}$, where $\beta\in\{3,4,8\}$.
Since in the coordinate ordering used above
none of the three codes has a generator matrix in standard form,
we apply suitable coordinate permutations.
The resulting generator matrices are 
$\left(
\begin{array}{cc}
I_{30} & M_\beta
\end{array}
\right)$, 
where the matrices
$M_\beta$ are listed in Table~\ref{Tab} and
$I_{30}$ denotes the identity matrix of order $30$.

\section{Concluding remarks}\label{Sec:R}

There is another reason for our interest in length $60$.
The codes constructed in this note belong to the family of binary
singly even self-dual $[10m,5m,2m]$ codes whose shadows have
minimum weight $m$, where $m$ is a positive integer.
We consider the weight enumerators that occur for codes in this family.

The general bound of Rains~\cite{Rains} implies that these parameters
can occur only for
\[
m\in\{1,2,3,4,5,6,7,8,10,12\}.
\]
For such a code, let $\beta$ denote the number of vectors of weight
$m$ in its shadow.
By~\cite[Lemmas~2 and~3]{HM}, we have
\[
\beta\in\{1,2,\ldots,8,10\}.
\]
Moreover, $\beta\le 2$ if $m$ is odd.
Here $\beta=0$ is excluded since the shadow is required
to have minimum weight $m$.

\begin{itemize}
\item
$m=1,2,3$:
The weight enumerators that occur were determined (see~\cite{C-S}).

\item
$m=4$:
The values of $\beta$ that occur are exactly
$1,2,\ldots,8$ and $10$~\cite{HM}.

\item
$m=5$:
The values of $\beta$ that occur are exactly $1$ and $2$
(see~\cite{DGH}, \cite{HM} and \cite{YL}).

\item
$m=6$:
Together with the previously known codes,
Proposition~\ref{main} shows that the values of $\beta$
that occur are exactly $1,2,\ldots,8$ and $10$.

\item
$m=7$:
If a binary self-dual $[70,35,14]$ code exists, then its shadow
has minimum weight $15$~\cite{DH}.
Hence, no such code has a shadow of minimum weight $7$.

\item
$m=8$:
It is not known whether a binary singly even self-dual
$[80,40,16]$ code exists.

\item
$m=10$:
There is no binary self-dual $[100,50,20]$ code~\cite{DGH}.

\item
$m=12$:
Every binary extremal self-dual $[120,60,24]$ code must be
doubly even~\cite{Rains}.
Hence, there is no binary singly even self-dual
$[120,60,24]$ code.
\end{itemize}

Consequently, the weight enumerators that occur in this family
have now been determined for every $m$ except $m=8$.
It remains to determine which weight enumerators occur for binary extremal
singly even self-dual $[80,40,16]$ codes whose shadows have
minimum weight $8$.
In fact, even the existence of a binary extremal singly even self-dual
$[80,40,16]$ code is not known.

\bigskip
\noindent
\textbf{Acknowledgments.}
This work was supported by JSPS KAKENHI Grant Number 23K25784.

\bigskip
\noindent
\textbf{Use of generative AI.}
The author used OpenAI's ChatGPT (GPT-6 Astra) to assist in constructing 
the codes in October 2026.
The author independently verified the computational
results presented in this note by \textsc{Magma}~\cite{Magma}.


\begin{landscape}

\begin{table}[thbp]
\caption{Matrices $M_{3}$, $M_{4}$ and $M_{8}$}
\label{Tab}
\begin{center}
{\footnotesize
\resizebox{\linewidth}{!}{%
\begin{tabular}{c}
$M_{3}=
\left(
\begin{array}{c}
101110001011011111111001110001 \\
101110010011110001111001101111 \\
101110111011000010111001000100 \\
110100011101111110010110110101 \\
100011011000100101001110101101 \\
111001011001001111100001111010 \\
001101100010101010000001011011 \\
010111010110110110000011010010 \\
011010001100011100000010110110 \\
101100010010001110000011010101 \\
101111111111111101000000111101 \\
101101001110100111000010110010 \\
011111110011000111101111000101 \\
011110011010000000110000111001 \\
011101101001000000011000111100 \\
000110000110110000000010011101 \\
010001000011011000000001101110 \\
001011000010101111000011001100 \\
111011110110011110000101110011 \\
110101101010110010000100010001 \\
011001111101100111010010110111 \\
001111111001010000001110100000 \\
010110111011001111100100101111 \\
011100011111011011101100010100 \\
011100101111101101110010001010 \\
011100110111110110011001100001 \\
000000001111110011110011011100 \\
000000100111011101011101101010 \\
000000010111101110101110110001 \\
000011000111000111000111111111 \\
\end{array}
\right)
$,
$M_{4}=
\left(
\begin{array}{c}
000111101111110100000111010011 \\
000111110111011010000111001101 \\
000111011111101001000111100110 \\
011101010110000100111110100001 \\
110001000011001001110000110000 \\
101011010101110010001110010001 \\
110001011101101010111110011011 \\
101011111110100111101101100100 \\
011101100011001010101011111111 \\
000101001001011101000101001000 \\
000110100100101110000110100000 \\
000011010010110011000011010000 \\
000011001011000000101111111101 \\
000101100101000000110111111110 \\
000110010110000000011111111011 \\
011101010110100001010011010011 \\
110001000011011111000110101001 \\
101011010101111001010101111101 \\
100110010101001010111011010001 \\
010011001110100001111101001100 \\
111010110011000010101000101011 \\
011110101010001001111100111001 \\
110100111001101111110001011000 \\
101101101100100001001010100001 \\
001000001000001010111101100010 \\
001000111000111100100011111100 \\
001000100000100111001000010111 \\
111111010000010110010110001001 \\
001000011000011101011101101101 \\
001000110000110011110011011011 \\
\end{array}
\right)
$,
$M_{8}=
\left(
\begin{array}{c}
100010000110101010111101111101 \\
110010100110001010010011100100 \\
111010110110011010100000001000 \\
101000010100110001000101011110 \\
101011010111011111101011110110 \\
101110010010101100011000000110 \\
000110110000001011110001000010 \\
000011011000100101011100000001 \\
000101101000010110101010000100 \\
010100011110010111101001011101 \\
001010101011001111110100101110 \\
011001001101011111011010110100 \\
010111010100111001011101000010 \\
001111001010111100101110000001 \\
011111011001000010110011000011 \\
101100001100011011001011010011 \\
101001111111101110001000010000 \\
101111100010110000001101010101 \\
110011100101111011101101011010 \\
111110110011000110011011110011 \\
010111000110000001000111110100 \\
001111000011000100000111011010 \\
011111111101111010000111101110 \\
011110101111100001000000010100 \\
010011001111101100000000001101 \\
001101100111110010000000100110 \\
011110000010001111010000110100 \\
010011111001011111001000011101 \\
001101111100101111100000101110 \\
011000111000111000000000000111 \\
\end{array}
\right)
$
\end{tabular}%
}
}
\end{center}
\end{table}

\end{landscape}

\end{document}